\documentclass[oneside]{amsart}

\usepackage[letterpaper,body={12.6cm,21cm}, mag=1000]{geometry}
\usepackage{amssymb}
\usepackage{amsthm}
\usepackage{amscd}

\numberwithin{equation}{section}

\theoremstyle{plain}

\newtheorem{thm}{Theorem}[section]
\newtheorem{cor}[thm]{Corollary}

\newtheorem{lemma}[thm]{Lemma}

\newtheorem{prop}[thm]{Proposition}
\newtheorem{step}{Step}
\newtheorem{remark}[thm]{Remark}

\theoremstyle{definition}

\newcommand{\dlabel}[1]{\ifmmode \text{\ttfamily \upshape [#1] } \else
{\ttfamily \upshape [#1] }\fi \label{#1}}

\newcommand{\B}{\operatorname{B} }
\newcommand{\C}{\operatorname{C} }

\newcommand{\Z}{\operatorname{Z} }
\newcommand{\A}{\operatorname{A} }

\newcommand{\gen}[1]{\left < #1 \right >}
\newcommand{\Aut}{\operatorname{Aut} }

\newcommand{\Hom}{\operatorname{Hom} }
\newcommand{\Inn}{\operatorname{Inn} }

\newcommand{\Out}{\operatorname{Out} }

\newcommand{\Autcent}{\operatorname{Autcent} }

\begin{document}

\title{On automorphisms of some finite $p$-groups}
\author{Manoj K.~Yadav}

\address{School of Mathematics, Harish-Chandra Research Institute \\
Chhatnag Road, Jhunsi, Allahabad - 211 019, INDIA}

\email{myadav@mri.ernet.in}

\thanks{Research supported by DST (SERC Division), the Govt. of INDIA}

\date{\today}

\begin{abstract}
 We give a sufficient condition on a finite 
$p$-group $G$ of nilpotency class $2$ so that $\Aut_c(G) = \Inn(G)$, where
$\Aut_c(G)$ and $\Inn(G)$ denote the group of all class preserving
automorphisms and inner automorphisms of $G$ respectively. Next we prove that
if $G$ and $H$ are two isoclinic finite groups (in the sense of P. Hall), then
$\Aut_c(G) \cong \Aut_c(H)$. Finally we study class preserving automorphisms
of groups of order $p^5$, $p$ an odd prime and prove that 
$\Aut_c(G) = \Inn(G)$ for all the groups $G$ of order $p^5$ except two 
isoclinism families.

\vspace{.2in}
{\bf Subj-class.} GR.  
\end{abstract}

\maketitle

\section{Introduction}

Let $G$ be a finite $p$-group and $|G| = p^n$, where $p$ is a prime and $n$ is
a positive integer. For $x \in G$, $x^G$ denotes the conjugacy class of 
$x$ in $G$. By $\Aut(G)$ we denote the group of all automorphisms of $G$.
An automorphism $\alpha$ of $G$ is called \emph{class preserving} if 
$\alpha(x) \in x^G$ for all $x \in G$. The set of all class preserving
automorphisms of $G$,  denoted by $\Aut_{c}(G)$, is a normal subgroup of 
$\Aut(G)$. Notice that $\Inn(G)$, the group of all inner automorphisms of $G$,
 is a normal subgroup of $\Aut_c(G)$. 
Let $\Out_c(G)$ denote the  group $\Aut_c(G)/\Inn(G)$.

In 1911, W. Burnside \cite[pg. 463]{wB55} posed the following question: 
Does there exist any  finite group $G$ such that $G$ has a non-inner 
class preserving automorphism?
In 1913,  Burnside \cite{wB13} himself gave an affirmative answer to this 
question. 
He constructed a group $G$ of order $p^6$, $p$ an odd prime, such that
$\Out_c(G) \not= 1$. In \cite{hH80}, \cite[pg. 102-103]{mH04}, \cite{JK75}, 
\cite{iM92}, \cite{cS68} and \cite{fS03}  more groups were
constructed such that $\Out_c(G) \not= 1$. But the order of all these groups
is $\ge p^6$. 
It follows from \cite{KV01} that $\Out_c(G) = 1$ for all the groups $G$ of 
order $p^4$. That there exist groups $G$ of order $2^5$ such that 
$\Out_c(G) \not= 1$ follows from \cite{KV02} or \cite{gW47}. In this paper we
study the class preserving automorphisms of groups of order $p^5$ for odd
primes and prove the following theorem:

\noindent{\bf Theorem} (Theorem \ref{thm3a}).
 Let $G$ be a finite $p$-group of order $p^5$, 
where $p$ is an odd prime. Then $\Out_c(G) \not= 1$ if and only if $G$ 
is isoclinic to  one of the groups $G_7$, $G_{10}$ and $H$, defined in 
\eqref{eqn1a}, \eqref{eqn1b} and \eqref{eqn1c} respectively.

Thus for any prime $p$, $n = 4$ is the largest number such that 
$\Out_c(G) = 1$ for all the groups $G$ of order $\leq$ $p^n$.

In section $4$ we prove that if $G$ and $H$ are two isoclinic finite groups 
(see Section $4$ below), then  $\Aut_c(G) \cong \Aut_c(H)$. 
This result allows us to study $\Aut_c(G)$ for a group $G$ only upto 
isoclinism. A list of the groups of order $p^5$ for odd primes $p$, ordered
into ten isoclinism families, is available from James' work \cite{rJ80}.
Our method of proof is to take one group $G$ from each isoclinism
family and compute the order of $\Aut_c(G)$ by using the upper and lower
bounds derived in Section $2$. In section $3$ we prove some results regarding
class preserving automorphisms of finite $p$-group of class $2$.

Our notation for objects associated with a finite multiplicative group $G$ is 
mostly standard.  We use $1$ to denote both the identity element of 
$G$ and the trivial subgroup $\{ 1 \}$ of $G$. The abelian group of all 
homomorphisms from an abelian group $H$ to an abelian group $K$ is denoted 
by $\Hom(H,K)$. We write $\gen x$ for the cyclic subgroup of $G$ generated 
by a given element  $x \in G$. 
To say that some $H$ is a subset or a subgroup of $G$ we write 
$H \subseteq G$ or $H \leq G$ respectively. To indicate, in addition, that 
$H$ is properly contained in $G$, we write $H \subset G$, $H < G$ 
 respectively.  If $x,y \in G$, then $x^y$ denotes the conjugate element 
$y^{-1}xy \in G$ and $[x,y] = [x,y]_G$ denotes the commutator 
$x^{-1}y^{-1}xy = x^{-1}x^y \in G$. If $x \in G$, then $x^G$ 
denotes the $G$-conjugacy class of all $x^w$, for $w \in G$, and $[x,G]$ 
denotes the set of all $[x,w]$, for $w \in G$. 
For $x \in G$, $\C_{H}(x)$ denotes the 
centralizer of $x$ in $H$, where $H \le G$. The center of $G$ will be 
denoted by $\Z(G)$. The Frattini subgroup of $G$ is denoted by $\Phi(G)$.

We write the subgroups in the lower central series of $G$ as $\gamma_n(G)$,
where $n$ runs over all strictly positive integers. 
And we write the subgroups in the upper central series of $G$ as $\Z_n(G)$,
where $n$ runs over all non-negative integers. 
We will be using the commutator identities
\[[x,yz] = [x,z][x,y]^z = [x,z][x,y][x,y,z]\]
and
\[[xy,z] = [x,z]^y[y,z] = [x,z][x,z,y][y,z]\]
without any reference.

\noindent{\bf Acknowledgements.} I thank the referee for his/her useful 
comments and suggestions. After having seen the results of this paper, 
he/she has suggested alternative proofs, which are given in Section 6, 
of Lemma \ref{lemma7}, Lemma \ref{lemma8}  and Lemma \ref{lemma9}. 

\section{Some Useful Lemmas}

An automorphism $\phi$ of a group $G$ is called   \emph{central} if 
$g^{-1}\phi(g) \in \Z(G)$ for all  $g \in G$. The set of all central 
automorphisms of $G$, denoted by $\Autcent(G)$, is a normal subgroup 
of $\Aut(G)$. It follows from \cite[Proposition 1.7]{cS68} that 
$\Autcent(G) = \C_{\Aut(G)}(\Inn(G))$.
Following \cite{AY65}, we shall say that a finite group $G$ is  
\emph{purely non-abelian} if it does not have a non-trivial abelian direct 
factor.

The following lemma follows from \cite{AY65}.
\begin{lemma}\label{lemma0}
Let $G$ be a purely non-abelian finite $p$-group. Then $|\Autcent(G)|$ 
$ = |\Hom(G/\gamma_2(G), \Z(G))|$.
\end{lemma}



\begin{lemma}\label{lemma2}
Let $G$ be a finite $p$-group such that $\Z(G) \subseteq [x,G]$ for all
$x \in G-\gamma_2(G)$. Then $|\Aut_c(G)| \ge |\Autcent(G)|
|G/\Z_2(G)|$.
\end{lemma}
\begin{proof}
Let $f \in \Autcent(G)$. Then $f$ fixes $\gamma_2(G)$ element wise. Therefore 
$f(x) = x \in x^G$ for all $x \in \gamma_2(G)$. So let
$x \in G- \gamma_2(G)$. Since $f \in \Autcent(G)$, $x^{-1}f(x) \in \Z(G)
\subseteq [x,G]$. Thus $x^{-1}f(x) = [x,g]$ for some $g \in G$. This implies
that $f(x) = g^{-1}xg \in x^G$. Thus $f(x) \in x^G$ for all $x \in G$, which
proves that $f \in \Aut_c(G)$. Thus $\Autcent(G) \le \Aut_c(G)$.
Since $\Inn(G) \le \Aut_c(G)$, it follows that
\[|\Aut_c(G)| \ge |\Autcent(G)| |\Inn(G)|/|\Autcent(G) \cap \Inn(G)|.\]
Set $\A = \Autcent(G) \cap \Inn(G)$. Then $A \cong \Z(\Inn(G))$.
And $\Z(\Inn(G)) \cong \Z(G/\Z(G))$ $ = \Z_2(G)/\Z(G)$.
Thus $\A \cong \Z_2(G)/\Z(G)$ and therefore $|\A| = |\Z_2(G)/\Z(G)|
=|\Z_2(G)|/|\Z(G)|$. Hence
\begin{eqnarray*}
|\Aut_c(G)| &\ge& |\Autcent(G)| |\Inn(G)|/|\Autcent(G) \cap \Inn(G)|\\ 
&=& (|\Autcent(G)| |G|/|\Z(G)|)/(|\Z_2(G)|/|\Z(G)|)\\
&=& |\Autcent(G)| |G|/|\Z_2(G)|.
\end{eqnarray*} 
This completes the proof of the  lemma. \hfill $\Box$

\end{proof}




\begin{lemma}\label{lemma3}
Let $G$ be a finite group. Let $\tau$ be an endomorphism of $G$ such that
$\tau (x) \in x^G$ for all $x \in G$. Then $\tau \in \Aut_c(G)$.
\end{lemma}
\begin{proof}
Such a $\tau$ obviously has trivial kernel. \hfill $\Box$

\end{proof}

The following lemma is \cite[Proposition 14.4]{mH04}.
\begin{lemma}\label{lemma4}
Let $G$ be a finite group and $H$ be an abelian normal subgroup of $G$ such 
that $G/H$ is cyclic. Then $\Out_c(G) = 1$.
\end{lemma}

\begin{lemma}\label{lemma4a}
Let $G$ be a finite $p$-group of order $p^n$ and $G$ have a cyclic subgroup
of  order $p^{n-2}$, where $p$ is an odd prime. Then $\Out_c(G) = 1$.
\end{lemma}
\begin{proof} 
If $G$ has a cyclic subgroup of order $p^{n-1}$, then the lemma follows from
Lemma \ref{lemma4}. Otherwise it follows from \cite{FN04} and \cite{KV02}.
\hfill $\Box$

\end{proof}

\begin{lemma}[\cite{mY06}]\label{lemma4b}
Let $G$ be a finite group. Let  $\{x_1, x_2, \dots, x_d\}$ be
a minimal generating set for $G$. Then $|\Aut_c(G)| \le \Pi_{i=1}^d |x_i^G|$.
\end{lemma}
\begin{proof}
There are no more than $\Pi_{i=1}^d |x_i^G|$ choices for the images of the
given generators to define a class preserving automorphism of $G$. 
\hfill $\Box$ 

\end{proof}

\section{Groups of class $2$}

Let $G$ be a finite nilpotent group of  class $2$. Let $\phi \in
\Aut_c(G)$. Then the map $g \mapsto g^{-1}\phi(g)$ is a homomorphism of
$G$ into $\gamma_2(G)$. This homomorphism sends $Z(G)$ to $1$. So it
induces a homomorphism $f_{\phi} \colon G/Z(G) \to \gamma_2(G)$, sending
$gZ(G)$ to $g^{-1}\phi(g)$, for any $g \in G$.  It is easily seen that
the map $\phi \mapsto f_{\phi}$ is a monomorphism of the group
$\Aut_c(G)$ into $\Hom(G/Z(G), \gamma_2(G))$.

Any $\phi \in \Aut_c(G)$ sends any $g \in G$ to some $\phi(g)
\in g^G$. Then $f_{\phi}(gZ(G)) = g^{-1}\phi(g)$ lies in $g^{-1}g^G =
[g,G]$.  Denote
\[  \{ \, f \in \Hom(G/Z(G), \gamma_2(G)) \mid f(gZ(G)) \in [g,G], \text{ for
all $g \in G$}\,\} \]
by $\Hom_c(G/Z(G), \gamma_2(G))$. Then $f_{\phi} \in \Hom_c(G/Z(G),
\gamma_2(G))$ for all $\phi \in \Aut_c(G)$. On the other hand, if $f \in
\Hom_c(G/Z(G), \gamma_2(G))$, then the map sending any $g \in G$ to
$gf(gZ(G))$ is an automorphism $\phi \in \Aut_c(G)$ such that $f_{\phi}
= f$. Thus we have

\begin{prop}\label{prop1}
 Let $G$ be a finite nilpotent group of class 2. Then the
above map $\phi \mapsto f_{\phi}$ is an isomorphism of the group
$\Aut_c(G)$ onto $\Hom_c(G/$ $Z(G), \gamma_2(G))$.
\end{prop}

The following lemmas are well known.
\begin{lemma}\label{lemma5}
Let $\A$, $\B$ and $\C$ be finite abelian groups. Then

{\em (i)} $\Hom(\A \times \B, \C) \cong \Hom(\A,\C) \times \Hom(\B,\C)$;

{\em (ii)}  $\Hom(\A,  \B \times \C) \cong \Hom(\A,\B) \times \Hom(\A,\C)$.
\end{lemma}

\begin{lemma}\label{lemma6}
Let $\C_n$ and $C_m$ be two cyclic groups of order $n$ and $m$
respectively. Then $\Hom(\C_n, \C_m) \cong \C_d$, where $d$ is the greatest
common divisor of $n$ and $m$, and $\C_d$ is the cyclic group of order $d$.
\end{lemma} 

Let $G$ be a finite $p$-group of class $2$. Notice that $[x,G]$ is a
non-trivial proper normal subgroup of $G$ for all $x \in G-\Z(G)$.
Let $\{x_1, x_2, \dots, x_d\}$ be a minimal generating set for $G$. 
Then $G/\Z(G) = \times_{i=1}^d
\gen{\bar{x_i}}$, where $\bar{x_i} = x_i \Z(G)$ and some of the factors may
possibly be trivial (this may happen in the case when $\Z(G) \not\le \Phi(G)$).
Let $f \in \Hom_c(G/Z(G), \gamma_2(G))$. So $f(g\Z(G)) \in
[g,G]$ for all $g \in G$. In particular, $f(x_i\Z(G)) \in
[x_i,G]$, $1 \le i \le d$. Thus it
follows that $|\Hom_c(G/\Z(G), \gamma_2(G))| \le 
\Pi_{i=1}^d |\Hom(\gen{\bar{x_i}}, [x_i,G])|$. Since there is an isomorphism
from $\Aut_c(G)$ onto $\Hom_c(G/Z(G),$ $ \gamma_2(G))$, we have the following.

\begin{prop}\label{prop2}
Let $G$ be  a finite $p$-group of class $2$ and $\{x_1, x_2, \dots, x_d\}$ be
a minimal generating set for $G$. 
Then $|\Aut_c(G)| \le \Pi_{i=1}^d |\Hom(\gen{\bar{x_i}}, [x_i,G])|$.
\end{prop}

\begin{thm}\label{thm2}
Let $G$ be a finite $p$-group of class $2$. Let $\{x_1, x_2, \dots, x_d\}$ be
a minimal generating set for $G$ such that $[x_i, G]$ is cyclic, $1 \le i \le
d$. Then $\Out_c(G) = 1$.
\end{thm}
\begin{proof}
Since $|\Hom(\gen{\bar{x_i}}, [x_i,G])| \le |\gen{\bar{x_i}}|$, it follows
from Proposition \ref{prop2} that 
$|\Aut_c(G)| \le \Pi_{i=1}^d |\gen{\bar{x_i}}|$. Now
\[|\Inn(G)| \le |\Aut_c(G)| \le \Pi_{i=1}^d |\gen{\bar{x_i}}| = |G/\Z(G)| =
|\Inn(G)|.\]
Thus $|\Aut_c(G)| = |\Inn(G)|$ and therefore $\Out_c(G) = 1$.
\hfill $\Box$

\end{proof}

\begin{cor}\label{cor1}
Let $G$ be a finite $p$-group of class $2$ such that $\gamma_2(G)$ is cyclic. 
Then $\Out_c(G) = 1$.
\end{cor}


\section{Isoclinic groups}

The following concept was introduced by P. Hall \cite{pH40} (also see
\cite[pg. 39-40]{DY06} for details).

Let $X$ be a finite group and $\bar{X} = X/\Z(X)$. 
Then commutation in $X$ gives a well defined map
$a_{X} : \bar{X} \times \bar{X} \mapsto \gamma_{2}(X)$ such that
$a_{X}(x\Z(X), y\Z(X)) = [x,y]$ for $(x,y) \in X \times X$.
Two finite groups $G$ and $H$ are called \emph{isoclinic} if 
there exists an  isomorphism $\phi$ of the factor group
$\bar G = G/\Z(G)$ onto $\bar{H} = H/\Z(H)$, and an isomorphism $\theta$ of
the subgroup $\gamma_{2}(G)$ onto  $\gamma_{2}(H)$
such that the following diagram is commutative
\[
 \begin{CD}
   \bar G \times \bar G  @>a_G>> \gamma_{2}(G)\\
   @V{\phi\times\phi}VV        @VV{\theta}V\\
   \bar H \times \bar H @>a_H>> \gamma_{2}(H).
  \end{CD}
\]
The resulting pair $(\phi, \theta)$ is called an \emph{isoclinism} of $G$ 
onto $H$. Notice that isoclinism is an equivalence relation among finite 
groups.

\begin{thm}\label{thm1}
Let $G$ and $H$ be two finite non-abelian isoclinic groups. Then
$\Aut_c(G) \cong \Aut_c(H)$.
\end{thm}
\begin{proof}
Since $G$ and $H$ are isoclinic, there exist isomorphisms $\phi : G/\Z(G) 
\mapsto H/\Z(H)$ and $\theta : \gamma_2(G) \mapsto \gamma_2(H)$ such that
$\theta([x,y]) = [x'\Z(H),y'\Z(H)] = [x',y']$, where $x'\Z(H) = \phi (x\Z(G))$
and  $y'\Z(H) = \phi (y\Z(G))$. 
 
Let $\tau \in \Aut_c(G)$. Let $x' \in H-\Z(H)$. Then $x'\Z(H)
 \in H/\Z(H)$ and $x\Z(G) = \phi^{-1}(x'\Z(H)) \in G/\Z(G)$. 
So $x \in G$ and therefore
there exists an element $a_x \in G$ such that $\tau (x) = a_x^{-1}xa_x$.

Let $a'_x$ be a coset representative of  $\phi(a_x\Z(G))$.
 Now
define a map $\sigma_{\tau} : H \mapsto H$ by
\begin{equation}
\label{eq1} \sigma_{\tau}(x') =
   \begin{cases}
     (a'_x)^{-1}x'a'_x      &\text{if $x' \in H-\Z(H)$ ;}\\
     x', &\text{if $x' \in \Z(H)$.}
   \end{cases}
\end{equation}
To make the proof more readable, we prove it in several steps.
\begin{step}
$\sigma_{\tau}$ is well defined.
\end{step}
\begin{proof}
We prove that for each $x' \in H$, $\sigma_{\tau}(x')$ is unique.
Let $w$ and $x$ be two coset representatives of $\phi^{-1}(x'\Z(H))$ such that
$w \not= x$. Then $w = xz$, where $z \in \Z(G)$. 
Now 
\[\tau(w) = \tau(xz) = \tau(x) \tau(z) = a_x^{-1}xa_x z = a_x^{-1}xza_x 
=a_x^{-1}wa_x.\]
Thus it follows that  $a_x$ in $\tau(x)$ is independent of the choice of coset
representative of the coset $x\Z(G)$.

 Now suppose that there exist two elements $a_x$ and $b_x$ in $G$
such that $\tau (x) = a_x^{-1}xa_x = b_x^{-1}xb_x$. 
Then, for $\sigma_{\tau}(x')$ there are two choices $(a'_x)^{-1}x'a'_x$ 
and $(b'_x)^{-1}x'b'_x$, where
$a'_x$ and $b'_x$ are coset representatives of $\phi(a_x\Z(G))$ and 
$\phi(b_x\Z(G))$ respectively. We claim that $(a'_x)^{-1}x'a'_x = 
(b'_x)^{-1}x'b'_x$. Since $a_x^{-1}xa_x = b_x^{-1}xb_x$, it follows that
$[x,b_xa_x^{-1}] = 1$. Now applying $\theta$ on it we get
\begin{eqnarray*}
1 &=& \theta([x,b_xa_x^{-1}]) = [\phi(x\Z(G)), 
\phi(b_xa_x^{-1}\Z(G))]\\
&=& [\phi(x\Z(G)), 
\phi(b_x\Z(G)) \phi(a_x^{-1}\Z(G))] = [x'\Z(H), b'_x\Z(H)(a'_x)^{-1}\Z(H)]\\
&=& [x'\Z(H), b'_x(a'_x)^{-1}\Z(H)] = [x', b'_x(a'_x)^{-1}].
\end{eqnarray*}
Thus it follows that $(a'_x)^{-1}x'a'_x = (b'_x)^{-1}x'b'_x$ and therefore our
claim is true.
Finally suppose that $w' \in a'_x\Z(H)$. 
Then $w' = a'_xz'$ for some $z' \in \Z(H)$. Now
\[(w')^{-1}x'w' = (a'_xz')^{-1}x'a'_xz' = (a'_x)^{-1}x'a'_x.\] 
Thus $\sigma_{\tau}(x')$ is independent of the choice of coset representative 
of $a'_x\Z(H)$. This proves that $\sigma_{\tau}$ is well defined.
\end{proof}

\begin{step}
$\sigma_{\tau} \in \Aut_c(H)$.
\end{step}
\begin{proof}
Let $x', y' \in H$. If both $x', y' \in \Z(H)$, then $\sigma_{\tau}(x'y') =
x'y' = \sigma_{\tau}(x') \sigma_{\tau}(y')$. Now let $x' \in H-\Z(H)$ and $y'
\in \Z(H)$. Then $x'y' \in H-\Z(H)$ and $x'y'\Z(H) = x'\Z(H)$. Thus 
$\sigma_{\tau}(x'y') = (a'_x)^{-1}x'y'a'_x = (a'_x)^{-1}x'a'_xy' =
\sigma_{\tau}(x')\sigma_{\tau}(y')$, since $y' \in \Z(H)$. So assume that 
$x', y' \in H- \Z(H)$.
Now 
\[\phi^{-1}(x'y'\Z(H)) = \phi^{-1}(x'\Z(H))
\phi^{-1}(y'\Z(H)) = x\Z(G) y\Z(G) = xy\Z(G),\]
where $x, y \in G$.
 Let  $\tau(xy) =
a_{xy}^{-1}xya_{xy}$, $\tau(x) = a_x^{-1}xa_x$ and  $\tau(y) = a_y^{-1}ya_y$.
Since $\tau(xy) = \tau (x) \tau(y)$, we get
\[a_{xy}^{-1}xya_{xy} = a_x^{-1}xa_x a_y^{-1}ya_y\]
or
\[(xy)^{-1}a_{xy}^{-1}xya_{xy} = y^{-1}x^{-1}a_x^{-1}xa_x 
(yy^{-1})a_y^{-1}ya_y\]
or
\[[xy,a_{xy}] = [x^y, (a_x)^y][y,a_y],\]
where $(a_x)^y = y^{-1}a_x y$.
Now applying $\theta$ on both the sides, we get
\[[x'y',a'_{xy}] = [(x')^{y'}, (a'_x)^{y'}][y',a'_y]\]
or
\[(a'_{xy})^{-1}x'y'a'_{xy} = (a'_x)^{-1}x'a'_x (a'_y)^{-1}y'a'_y.\]
Thus from the definition \eqref{eq1} it follows that $\sigma_{\tau}(x'y') =
\sigma_{\tau}(x')\sigma_{\tau}(y')$. 
Hence $\sigma_{\tau}$ is an endomorphism of $H$. 
That $\sigma_{\tau} \in \Aut_c(H)$ follows from Lemma \ref{lemma3}. 
This completes the proof of Step $2$.
\end{proof}

\begin{step}
The map $\tau \mapsto \sigma_{\tau}$ is a homomorphism from $\Aut_c(G)$ to
$\Aut_c(H)$. 
\end{step}
\begin{proof}
Let $\tau_1, \tau_2 \in \Aut_c(G)$. Let $x' \in H$. Let $x$ be  a coset
representative of $\phi^{-1}(x'\Z(H))$. Since $\tau_1\tau_2 \in \Aut_c(G)$,
there exists $a_x \in G$ such that $\tau_1\tau_2(x) = a_x^{-1}xa_x$. Also
there exists $b_x, c_x \in G$ such that $\tau_2(x) = b_x^{-1}xb_x$ and 
$\tau_1(\tau_2(x)) = c_x^{-1}(b_x^{-1}xb_x)c_x$. Since $\tau_1\tau_2(x) = 
\tau_1(\tau_2(x))$, we get
\[a_x^{-1}xa_x = c_x^{-1}b_x^{-1}xb_xc_x\]
or
\[x^{-1}a_x^{-1}xa_x = x^{-1}(b_xc_x)^{-1}xb_xc_x\]
or
\[[x,a_x] = [x, b_xc_x].\]
Now applying $\theta$ on both the sides, we get
\[[x', a'_x] = [x', b'_xc'_x],\]
since $\phi(b_xc_x\Z(G)) = \phi(b_x\Z(G))\phi(c_x\Z(G)) = b'_x\Z(H)c'_x\Z(H)
= b'_x c'_x\Z(H)$,
where $a'_x$, $b'_x$ and $c'_x$ are the coset representatives of
$\phi(a_x\Z(G))$, $\phi(b_x\Z(G))$ and $\phi(c_x\Z(G))$ respectively.
Thus from the last equality we have 
\[(a'_x)^{-1}x'a'_x = (c'_x)^{-1}(b'_x)^{-1}x'b'_xc'_x.\]
Now it follows from the definitions of $\sigma_{\tau_1\tau_2}$,
  $\sigma_{\tau_1}$ and  $\sigma_{\tau_2}$ that 
$\sigma_{\tau_1\tau_2} (x')= (a'_x)^{-1}x'a'_x$ and 
$\sigma_{\tau_1}(\sigma_{\tau_2}(x')= (c'_x)^{-1}(b'_x)^{-1}x'b'_xc'_x$, for
all $x' \in H$. 
Hence $\sigma_{\tau_1\tau_2} = \sigma_{\tau_1}\sigma_{\tau_2}$. 
This proves  Step $3$.
\end{proof}

Similarly for each $\sigma \in \Aut_c(H)$ we can define $\tau_{\sigma} \in
\Aut_c(G)$ and the map sending $\sigma$ to $\tau_{\sigma}$ is a homomorphism
from $\Aut_c(H)$ to $\Aut_c(G)$. It is not difficult to prove that for each
$\tau \in \Aut_c(G)$ and $\sigma \in \Aut_c(H)$, $\tau_{\sigma_{\tau}} = \tau$
and $\sigma_{\tau_{\sigma}} = \sigma$. Thus it follows that the homomorphism 
from $\Aut_c(G)$ to $\Aut_c(H)$, defined above, becomes an isomorphism. 
This completes the proof of the
theorem. \hfill $\Box$

\end{proof}

\section{Groups of order $p^5$}

We'll use the classification of groups of order $p^5$ by R. James 
\cite[Section 4.5]{rJ80}.
Throughout this section $p$ always denotes an odd prime.

\begin{lemma}\label{lemma7}
Let $G$ be the group $\phi_7(1^5)$ in the isoclinism family (7) of
\cite[Section 4.5]{rJ80}. Then $|\Aut_c(G)| = p^5$.
\end{lemma}
\begin{proof}
Let $G$ be the group $\phi_7(1^5)$. Then $G$ is 
a nilpotent group of  class $3$ such that $\Z(G) \le \gamma_2(G) = \Phi(G)$, 
$|\gamma_2(G)| = p^2$ and $|\Z(G)| = p$.
Now it follows from \cite[Theorem 4.7]{TM04} that  $\Z(G) \le
[x,G]$ for all $x \in G-\Z(G)$. 
Since $\Z(G) \le \gamma_2(G)$, therefore $\Z(G) \le [x,G]$ for all 
$x \in G - \gamma_2(G)$. 
Thus it follows from Lemma \ref{lemma2} that
\begin{equation}\label{eqn1}
|\Aut_c(G)| \ge |\Autcent(G)| |G/\Z_2(G)| =
 \frac{|\Autcent(G)|p^5}{|\Z_2(G)|}.
\end{equation} 
From Lemma \ref{lemma0} we have $|\Autcent(G)| =
|\Hom(G/\gamma_2(G), \Z(G))| = p^3$, since $G/\gamma_2(G)$ is an elementary
abelian group of order $p^3$. 
It follows that $|\Z_2(G)| = p^3$. 
Thus from \eqref{eqn1} we
get $|\Aut_c(G)| \ge p^5$. It follows from \cite[Section 4.1]{rJ80} that 
in $G$ there are $p^3 -p$ conjugacy classes of length $p^2$. Thus all these
conjugacy classes covers $p^5-p^3$ elements out of $p^5-p^2$ elements in
$G-\gamma_2(G)$. So there must exists an element $y \in G-\gamma_2(G)$ such
that $|y^G| \le p$. Since $\Z(G) \le \gamma_2(G)$, $y$ is not a central
element and therefore $|y^G| = p$. We can extend the set $\{y\}$ to get a 
minimal generating set (say)  $\{y, x_1, x_2\}$ of $G$, 
since $\gamma_2(G) = \Phi(G)$.
Then from Lemma \ref{lemma4b} we get $|\Aut_c(G)| \le p^5$, since $|x_i^G| \le
p^2$, $i = 1,\;2$. Hence  $|\Aut_c(G)| = p^5$.
\hfill $\Box$

\end{proof}

\begin{lemma}\label{lemma8}
Let $G$ be the group $\phi_{10}(1^5)$ in the isoclinism family (10) of
\cite[Section 4.5]{rJ80}. Then $|\Aut_c(G)| = p^5$.
\end{lemma}
\begin{proof}
The group $G$ is a $p$-group of maximal class. $G$ is generated by $\alpha$
and $\alpha_1$ such that  the elements $\alpha_{i+1} := [\alpha_i, \alpha]$,
$1 \le i \le 3$, generate $\gamma_2(G)$ and $[\alpha_1, \alpha_2] = \alpha_4$.
Here $\alpha_2$,  $\alpha_3$ and $\alpha_4$ commutes with one another. Thus
$\gamma_2(G)$ is abelian. Also $[\alpha_1, \alpha_3] = 1$.
It is easy to prove that every element $g \in G$ can be written as
\[g = \beta \alpha_1^{k_1}\alpha^{k},\]
where $0 \le k_1, k \leq p-1$ and  $\beta \in \gamma_2(G)$.
Let $k \not= 0$. Then 
\[[g, \alpha_3] = [\beta \alpha_1^{k_1}\alpha^{k},\alpha_3] = [\alpha^k,
  \alpha_3] = \alpha_4^{-k} \in \Z(G),\]
since $[\beta \alpha_1^{k_1}, \alpha_3] = 1$.
Now let $k =0$ and $k_1 \not= 0$. Then
\[[g, \alpha_2] = [\beta \alpha_1^{k_1}, \alpha_2] = [\alpha_1^{k_1},
  \alpha_2] = \alpha_4^{k_1} \in \Z(G).\]
Thus it follows that $\Z(G) \subseteq [g,G]$ for $g \in G-\gamma_2(G)$,
since $|\Z(G)| = p$. Now from Lemma \ref{lemma2} we have
\[|\Aut_c(G)| \ge \frac{|\Autcent(G)|p^5}{|\Z_2(G)|} = p^5,\]
since $|\Z_2(G)| = p^2$ and $|\Autcent(G)| = |\Hom(G/\gamma_2(G), \Z(G))| =
p^2$. Since there are only $p^2-p$ conjugacy classes of length $p^3$, there
must exists an element $y \in G-\gamma_2(G)$ such that $|y^G| \le p^2$. We can
always extend $\{y\}$ to some minimal generating set (say) $\{y, x\}$ of $G$, 
therefore we have from Lemma \ref{lemma4b} that $|\Aut_c(G)| \le |y^G||x^G|
 \le p^5$. Hence $|\Aut_c(G)| = p^5$.
\hfill $\Box$

\end{proof}

\begin{lemma}\label{lemma9}
Let $G$ be the group $\phi_{6}(1^5)$ in the isoclinism family (6) of
\cite[Section 4.5]{rJ80}. Then $\Out_c(G) = 1$.
\end{lemma}
\begin{proof}
$G$ is a group of class $3$ and $|\Z(G)| = p^2$.
Suppose that $\Out_c(G) \not= 1$. Since $|\Inn(G)| = |G/\Z(G)| = p^3$, we must
have $|\Aut_c(G)| \ge p^4$. It follows from \cite[Section 4.1]{rJ80}
that $|x^G| = p^2$ for all $x \in
G-\Z(G)$. Also $|\gamma_2(G)| = |\Phi(G)| = p^3$. Then $G$ can be generated by
two elements (say) $x, y$ of $G$. Now it follows from Lemma \ref{lemma4b} that
$|\Aut_c(G)| \le |x^G||y^G| \le p^4$. Since $|\Inn(G)| = p^3$ and 
 $\Out_c(G) \not= 1$, $|\Aut_c(G)|$ must be $p^4$. Thus it follows that for
any two elements $u \in x^G$ and $v \in y^G$, there must exist an  
automorphism $\tau \in \Aut_c(G)$ such that $\tau(x) = u$ and $\tau(y) = v$.

The group $G$ is generated by $\alpha_1$, $\alpha_2$ such that 
$\gamma_2(G)$ is generated by $\beta := [\alpha_1,\alpha_2]$, $\beta_1 :=
[\beta, \alpha_1]$ and $\beta_2 := [\beta, \alpha_2]$, and 
$\Z(G)$ is generated by $\beta_1$ and $\beta_2$. 
Since $[\alpha_1, \alpha_2] \not= 1$, 
$\alpha_1 \not= \alpha_2^{-1}\alpha_1\alpha_2 \in \alpha_1^G$ and 
$\alpha_2 \not= \alpha_1^{-1}\alpha_2\alpha_1 \in \alpha_2^G$. 
Let $\tau \in \Aut_c(G)$ such
that $\tau(\alpha_1) = \alpha_2^{-1}\alpha_1\alpha_2$ and 
 $\tau(\alpha_2) = \alpha_1^{-1}\alpha_2\alpha_1$. Since  $\tau \in
\Aut_c(G)$, $(\alpha_1 \alpha_2)^{-1}\tau(\alpha_1 \alpha_2) \in 
[\alpha_1 \alpha_2, G]$. Now
\[(\alpha_1 \alpha_2)^{-1}\tau(\alpha_1 \alpha_2) = \alpha_2^{-1}\alpha_1^{-1}
\alpha_2^{-1}\alpha_1\alpha_2 \alpha_1^{-1}\alpha_2\alpha_1 = \beta_2.\]
Thus $\beta_2 \in [\alpha_1 \alpha_2, G]$. So there is an element $y \in G$
such that $\beta_2 = [\alpha_1 \alpha_2, y]$.
 We claim that $[\alpha_1 \alpha_2, \beta] \not= 1$. 
For, if $[\alpha_1 \alpha_2, \beta] =
1$, then $\gamma_2(G) < \C_{G}(\alpha_1 \alpha_2)$, since 
$\alpha_1 \alpha_2 \in G -\gamma_2(G)$. 
Thus $|\C_{G}(\alpha_1 \alpha_2)| \ge p^4$, which is a
contradiction to the fact that $|[\alpha_1 \alpha_2, G]| = p^2$. 
Hence our claim is true. Now $[\alpha_1 \alpha_2, \beta] 
= [\alpha_1, \beta][\alpha_2, \beta] = \beta_1^{-1}\beta_2^{-1}$, 
since $[\alpha_i, \beta] \in \Z(G)$, $i = 1, 2$. 
Thus $\beta_1^{-1}\beta_2^{-1} \in [\alpha_1 \alpha_2, G]$. Now
\[\beta_1^{-1} = \beta_1^{-1}\beta_2^{-1}\beta_2 = [\alpha_1 \alpha_2, \beta]
[\alpha_1 \alpha_2, y] = [\alpha_1 \alpha_2, \beta y].\]
Thus it follows that $\beta_1, \beta_2 \in [\alpha_1 \alpha_2, G]$ and 
therefore $\Z(G) \subseteq [\alpha_1 \alpha_2, G]$. 
Since $|[\alpha_1 \alpha_2, G]| = p^2$,
we get $\Z(G) = [\alpha_1 \alpha_2, G]$. This gives a contradiction, since 
$[\alpha_1 \alpha_2, \alpha_1] = [\alpha_2, \alpha_1] = \beta^{-1} \in
[\alpha_1 \alpha_2, G]$, but $\beta^{-1} \not\in \Z(G)$. This completes the
proof of the lemma.
\hfill $\Box$ 

\end{proof}

Define a set of relations $\mathcal R$ by
\begin{eqnarray*}
{\mathcal R} &=& \{a^p = b^p = x^p = y^p = z^p = 1\} \cup \{[x,y] = [x,z] = 
                 [y,z] = 1\}\\
          & & \quad \cup \{x^b = xz, y^b = y, z^b = z\} \cup \{x^a = xy, y^a =
                 yz, z^a = z\}.
\end{eqnarray*}

Now set 
\begin{equation}\label{eqn1a}
G_7 = \gen{a, b, x, y, z \;|\; {\mathcal R}, [b,a] = 1} \;\text{for}\; p \ge 3
\end{equation}
and 
\begin{equation}\label{eqn1b}
G_{10} = \gen{a, b, x, y, z \;|\; {\mathcal R}, [b,a] = x} \;\text{for}\; 
p \ge 5.
\end{equation}
The groups $G_7$ and $G_{10}$ are of the form $\gen{a} \ltimes (\gen{b} \ltimes
(\gen{x} \times \gen{y} \times \gen{z}))$ and have order $p^5$. 
For $p = 3$ define a group $H$ of order $3^5$ by
\begin{eqnarray}\label{eqn1c}
H &=& \gen{a, b, c \;| \;a^3 = b^9 = c^9 = 1, [b,c] = c^3, [a,c] = b^3, [b, a] 
= c}\\
&=& \gen{a} \ltimes (\gen{b} \ltimes \gen{c}).\nonumber
\end{eqnarray}

\begin{remark}\label{rem2}
The group  $\phi_7(1^5)$ in the isoclinism family (7) of
\cite[Section 4.5]{rJ80} is isomorphic to $G_7$. The group $\phi_{10}(1^5)$ 
in the isoclinism family (10) of \cite[Section 4.5]{rJ80} is isomorphic to
$G_{10}$ for $p \ge 5$ and is isomorphic to $H$ for $p = 3$.
\end{remark}

Now we prove our main theorem.
\begin{thm}\label{thm3a}
Let $G$ be a finite $p$-group of order $p^5$, where $p$ is an odd prime. Then
$\Out_c(G) \not= 1$ if and only if $G$ is isoclinic to one of the groups 
$G_7$, $G_{10}$ and $H$.
\end{thm}
\begin{proof}In view of Theorem \ref{thm1}, it is sufficient to study
  $\Out_c(G)$ only for one group $G$ from each isoclinism family of groups of
  order $p^5$. If $G$ is abelian, then obviously $\Out_c(G) = 1$. 
 Let $G$ be either $\phi_2(311)a$, 
 $\phi_3(311)a$ or $\phi_8(32)$. Then there exists an element $x \in G$
such that order of $x$ is $p^3$. 
So it follows from Lemma \ref{lemma4a} that $\Out_c(G) = 1$. 
Now let $G$ be the group $\phi_4(1^5)$. Then $G$ has a maximal abelian
subgroup $H$ such that $G/H$ is cyclic. Thus it follows from Lemma
\ref{lemma4} that $\Out_c(G) = 1$. If $G$ is some group from the
isoclinism family (5), then the class of $G$ is $2$ and $\gamma_2(G)$ is
cyclic. Thus from Corollary \ref{cor1} we have that $\Out_c(G) = 1$. 
Next consider any group $G$ from the isoclinism family (9). Then it follows
from \cite[Section 4.1]{rJ80} that the
nilpotency class of $G$ is $4$ and it has $p^2-p$ conjugacy classes of 
length $p^3$, $p^3-1$ conjugacy classes of length $p$ and $p$ conjugacy 
classes of length $1$. Thus $|\gamma_2(G)| = |\Phi(G)| = p^3$.
Any conjugacy class of length $p^3$ must be contained in
$G-\gamma_2(G)$, since $|\gamma_2(G)| = p^3$. Let $\C = \{x \in G | |x^G| =
p^3\}$. Then $|\C| = p^3 (p^2 - p) = p^5 -p^4$ and $\C \subset G-\gamma_2(G)$.
Since $|G-\gamma_2(G)| = p^5-p^3$, there must exist an element 
$y \in G-\gamma_2(G)$ such that $|y^G| < p^3$. Thus $|y^G| = p$. Since
$\gamma_2(G) = \Phi(G)$, $y \in  G-\Phi(G)$. Therefore the set $\{y\}$ can be
extended to a minimal generating set (say) $\{y,x\}$ of $G$. Then $|y^G||x^G|
\le p^4$. Now it follows from Lemma \ref{lemma4b} that 
$\Aut_c(G) \le |y^G||x^G|$. Thus
\[\Inn(G) \le \Aut_c(G) \le |y^G||x^G| \le  p^4 = |G/\Z(G)| = |\Inn(G)|.\]
Hence $\Inn(G) = \Aut_c(G)$, which gives $\Out_c(G) = 1$. 
Now it follows, along with Lemma \ref{lemma9}, that if $\Out_c(G) \not= 1$,
then $G$ can not lie in the isoclinism families (1) - (6), (8) and (9). 
Thus $G$ must lie either in  the isoclinism family (7) or (10). 
Hence $G$ is isoclinic to $\phi_7(1^5)$ or $\phi_{10}(1^5)$. Thus it follows
from Remark \ref{rem2} that $G$ is isoclinic to $G_7$, $G_{10}$ or $H$. 

Conversely suppose that $G$ is isoclinic to $G_7$, $G_{10}$ or $H$. Thus $G$
is isoclinic to $\phi_7(1^5)$ or $\phi_{10}(1^5)$.
Then it follows from Lemma \ref{lemma7} and Lemma \ref{lemma8} that $\Out_c(G)
\not= 1$. 
\hfill $\Box$

\end{proof}

\section{Some alternative proofs}

The following alternative proofs of some of our lemmas were provided by the
referee.

\begin{proof}[Alternative proof of Lemma \ref{lemma7} and Lemma \ref{lemma8}] 
Let $G$ be one of the groups $G_7$  and $G_{10}$.  Then the center of $G$ is
generated by $z$ and has order $p$. Let $\tau \in \Aut_c(G)$.
 Since $\Aut_c(G) = \Inn(G)$ for all groups $G$ of order $p^4$ (\cite{KV01}), 
$\tau$ induces an inner automorphism on $G/\Z(G)$ given by conjugation with 
(say) $a\Z(G)$. Thus $\tau(x) \Z(G) = x^a\Z(G)$ for all $x \in G$. 
So for each $x \in G$, there
 exists some element $z_x \in \Z(G)$ such that $\tau(x) = x^az_x = (xz_x)^a$.
Let us define a map $\tau'$ from $G$ to $G$ such that $\tau'(x) = a
\tau(x)a^{-1}$ for all $x \in G$. Then it is fairly easy to prove that 
$\tau' \in \Aut_c(G)$ and $\tau'(x) = xz_x$, $x \in G$. Now it follows that
$\tau'$ is a central automorphism of $G$ and $\tau = i_a\tau'$, where $i_a$
denotes the inner automorphism given by conjugation with $a$. Thus we have  
 $\Aut_c(G) \le \Autcent(G) \Inn(G)$. 
It follows (as in \cite[Section 3]{TM04}) that
$\Z(G) \le [g,G]$ for all $g \in G - \Z(G)$. Thus we have $\Aut_c(G) =
\Autcent(G) \Inn(G)$. Now using Lemma \ref{lemma0} one can prove that
$|\Autcent(G) \Inn(G)| = p^5$. Hence $|\Aut_c(G)| = p^5$. It is not difficult
to prove the result for the group $H$.   \hfill $\Box$

\end{proof}

\begin{proof}[Alternative proof of Lemma \ref{lemma9}]
The nilpotency class of $G$ is $3$ and $G$ has  an elementary abelian center
of order $p^2$. Let $\tau \in \Aut_c(G)$. Let $N_1, \cdots, N_{p+1}$ be the
central subgroups of order $p$ in $G$.
Since $\Aut_c(G) = \Inn(G)$ for all groups $G$ of order $p^4$ (\cite{KV01}), 
$\tau$ induces an inner automorphism on
$G/N_i$, $1 \le i \le p+1$.  Let $x_1, \cdots, x_{p+1} \in G$ such that 
$\tau$ induces on $G/N_i$ the inner automorphism given by
conjugation with $x_iN_i$. These inner automorphisms agree on $G/\Z(G)$, so
the $p+1$ elements $x_1^{-1}x_i$ all lie in $\Z_2(G)$. Since $|\Z_2(G)/\Z(G)|
= p$, we have $x_1^{-1}x_s\Z(G) = x_1^{-1}x_t\Z(G)$ for distinct indices $s$
and $t$. But then  $x_s\Z(G) = x_t\Z(G)$ and $\tau$ is the inner automorphism
given by conjugation with $x_s$. This proves that $\Out_c(G) = 1$. 
\hfill $\Box$ 

\end{proof}

\end{document}